\documentclass{article}
\usepackage{amssymb,latexsym,amsmath,amsthm,mathrsfs}
\usepackage{graphicx}
\newcommand{\comment}[1]{}

\begin{document}
\title{Various analytic observations on combinations\footnote{Presented to the
St. Petersburg Academy on April 6, 1741.
Originally published as
{\em Observationes analyticae variae de combinationibus},
Commentarii academiae scientiarum Petropolitanae \textbf{13} (1751), 64--93.
E158 in the Enestr{\"o}m index.
Translated from the Latin by Jordan Bell,
Department of Mathematics, University of Toronto, Toronto, Canada.
Email: jordan.bell@utoronto.ca}}
\author{Leonhard Euler}
\date{}
\maketitle

1. Let a series, either finite or infinite, of quantities be given to us,
as
\[
a,b,c,d,e,f,g,h \quad \textrm{etc.};
\]
the letters
denote these quantities, which may be either equal or unequal to each other.
At the same time, however, I will speak of quantities indicated by different
letters as being unequal to each other, even if in examples equal numbers
can be substituted in place of them.

2. Now first let new series be formed
from these quantities
by taking powers, with the sums designated by capital letters
$A,B,C,D$ etc. as follows; thus:
\begin{eqnarray*}
A&=&a+b+c+d+e+\textrm{etc.},\\
B&=&a^2+b^2+c^2+d^2+e^2+\textrm{etc.},\\
C&=&a^3+b^3+c^3+d^3+e^3+\textrm{etc.},\\
D&=&a^4+b^4+c^4+d^4+e^4+\textrm{etc.},\\
E&=&a^5+b^5+c^5+d^5+e^5+\textrm{etc.}\\
&&\textrm{etc.}
\end{eqnarray*}
These series will be infinite if the number of quantities $a,b,c,d$ etc.
has been assumed to be infinite; on the other hand if the number of these
quantities is finite and determinate, put $=n$, then all these series
will be comprised of that same number of terms.

3. Next  let us now form series by taking products
of unequal terms from the assumed quantities $a,b,c,d$ etc. in
the following way.
Namely, first the single quantities are summed,
then the products of two unequal terms, from three unequal terms, from
four unequal terms, and so on; and let us indicate these series with the Greek
letters $\alpha,\beta,\gamma,\delta$ etc., so that it follows:
\begin{eqnarray*}
\alpha&=&a+b+c+d+\textrm{etc.},\\
\beta&=&ab+ac+ad+ae+bd+\textrm{etc.},\\
\gamma&=&abc+abd+abe+bcd+\textrm{etc.},\\
\delta&=&abcd+abce+bcde+\textrm{etc.},\\
\epsilon&=&abcde+\textrm{etc.}\\
&&\textrm{etc.}
\end{eqnarray*}
If the number of assumed quantities $a,b,c,d$ etc. were infinite, then
these series would not only all extend infinitely, but also, the number
of forms of these series would be infinite.
But if instead the number
of the quantities $a,b,c,d$ etc. is finite, put $=n$, then
the series $\alpha$ will contain $n$ terms, the second series $\beta$ will
be comprised from $\frac{n(n-1)}{1\cdot 2}$ terms, the third
$\gamma$ from $\frac{n(n-1)(n-2)}{1\cdot 2\cdot 3}$ terms, the fourth
$\delta$ from $\frac{n(n-1)(n-2)(n-3)}{1\cdot 2\cdot 3\cdot 4}$ terms,
and so on, until finally a series is reached
that is comprised from a single term, and then 
all the subsequent series would vanish, as they could not have any terms.
It is also clear that the number of series that occur here is $=n$, the last
of which consists of a single term, which is the product of all
the assumed quantities $a,b,c,d,e$ etc.

4. Since here we only took products from unequal quantities and
formed the series treated above from them,
by thus repeating the same quantities in products,
new series of products from one, two, three, four etc. will be obtained,
in which equal factors are not excluded as before;
so the series will be obtained thus: 
\begin{eqnarray*}
\mathfrak{A}&=&a+b+c+d+e+\textrm{etc.},\\
\mathfrak{B}&=&a^2+ab+b^2+ac+bc+c^2+\textrm{etc.},\\
\mathfrak{C}&=&a^3+a^2b+ab^2+b^3+a^2c+abc+\textrm{etc.},\\
\mathfrak{D}&=&a^4+a^3b+a^2b^2+a^2bc+abcd+\textrm{etc.},\\
\mathfrak{E}&=&a^5+a^4b+a^3b^2+a^3bc+a^2bcd+\textrm{etc.}\\
&&\textrm{etc.};
\end{eqnarray*}
in other words, these series contain all the quantities 
which can be produced from the multiplication of the given quantities
$a,b,c,d$ etc. Again we should note that if the number of quantities
$a,b,c,d$ etc. were finite, $=n$, then the first series
$\mathfrak{A}$ must have $n$ terms; while the second $\mathfrak{B}$ would have
$\frac{n(n+1)}{1\cdot 2}$ terms, the third $\mathfrak{C}$ would have
$\frac{n(n+1)(n+2)}{1\cdot 2\cdot 3}$ terms, the fourth $\mathfrak{D}$
indeed $\frac{n(n+1)(n+2)(n+3)}{1\cdot 2\cdot 3\cdot 4}$ terms, and so on.

5. The three orders of series, which we have composed from the given
quantities $a,b,c,d$ etc. in three  ways, are connected to each other,
such  that with the order of any one of the series known, then from it 
the orders of the remaining two can be determined. For the details
of how this law works
and a method for investigating
this, observation and induction are typically applied for the most part;
first noticing of course that $A=\alpha=\mathfrak{A}$,
and for the remaining it has been checked that:
\begin{eqnarray*}
\alpha&=&A,\\
\beta&=&\frac{\alpha A-B}{2},\\
\gamma&=&\frac{\beta A-\alpha B+C}{3},\\
\delta&=&\frac{\gamma A-\beta B+\alpha C-D}{4},\\
\epsilon&=&\frac{\delta A-\gamma B+\beta C-\alpha D+E}{5}\\
&&\textrm{etc.}
\end{eqnarray*}
similarly
\begin{eqnarray*}
\mathfrak{A}&=&A,\\
\mathfrak{B}&=&\frac{\mathfrak{A}A+B}{2},\\
\mathfrak{C}&=&\frac{\mathfrak{B}A+\mathfrak{A}B+C}{3},\\
\mathfrak{D}&=&\frac{\mathfrak{C}A+\mathfrak{B}B+\mathfrak{A}C+D}{4},\\
\mathfrak{E}&=&\frac{\mathfrak{D}A+\mathfrak{C}B+\mathfrak{B}C+\mathfrak{A}D+E}{5}\\
&&\textrm{etc.}
\end{eqnarray*}
and also
\begin{eqnarray*}
\mathfrak{A}&=&\alpha,\\
\mathfrak{B}&=&\alpha \mathfrak{A}-\beta,\\
\mathfrak{C}&=&\alpha \mathfrak{B}-\beta\mathfrak{A}+\gamma,\\
\mathfrak{D}&=&\alpha \mathfrak{C}-\beta\mathfrak{B}+\gamma\mathfrak{A}-\delta,\\
\mathfrak{E}&=&\alpha\mathfrak{D}-\beta\mathfrak{C}+\gamma\mathfrak{B}-\delta\mathfrak{A}+\epsilon\\
&&\textrm{etc.}
\end{eqnarray*}
By means of these relations, given the the sums of the series in any one
class, the sums of the series contained in the other two classes will be able
to be defined.

6. By carefully studying the nature and qualities of these series,
the truth of this mutual relation will indeed easily be clear by observation
and induction. Truly though, however much we are convinced by 
the truth of this connection, it will be fruitful to consider
the entire problem in the following way;
whence at once other properties are offered to us in addition
which induction alone
does not easily
present a path to. Namely, assuming as
given the quantities
\[
a,b,c,d,e \quad \textrm{etc.}
\]
from which are formed the three classes of series detailed above,
let us consider this expression
\[
P=\frac{az}{1-az}+\frac{bz}{1-bz}+\frac{cz}{1-cz}+\frac{dz}{1-dz}+\frac{ez}{1-ez}+\textrm{etc.},
\]
and with all the terms resolved into geometric progressions in the usual
way this gives
\begin{eqnarray*}
P&=&+z(a+b+c+d+e+\textrm{etc.})\\
&&+z^2(a^2+b^2+c^2+d^2+e^2+\textrm{etc.})\\
&&+z^3(a^3+b^3+c^3+d^3+e^3+\textrm{etc.})\\
&&+z^4(a^4+b^4+c^4+d^4+e^4+\textrm{etc.})\\
&&\textrm{etc.};
\end{eqnarray*}
all these series are contained in the first class.
Then if the sums given above (\S 2) are written in place
of these, it will be
\[
P=Az+Bz^2+Cz^3+Dz^4+Ez^5+\textrm{etc.},
\]
and thus the sum of this series will be, as we assumed,
\[
P=\frac{az}{1-az}+\frac{bz}{1-bz}+\frac{cz}{1-cz}+\frac{dz}{1-dz}+\textrm{etc.}
\]

Also in a similar way, if
\[
Q=\frac{az}{1+az}+\frac{bz}{1+bz}+\frac{cz}{1+cz}+\frac{dz}{1+dz}+\textrm{etc.},
\]
by series of the first class it will be
\[
Q=Az-Bz^2+Cz^3-Dz^4+Ez^5-\textrm{etc.}
\]

7. Let us next consider this expression
\[
R=(1+az)(1+bz)(1+cz)(1+dz)(1+ez)\, \textrm{etc.};
\]
if the factors are actually multiplied into each other and
the terms are disposed according to the exponents of $z$,
the coefficient of $z$ will be equal to the
sum of the given quantities $a,b,c,d$ etc.
The coefficient of $z^2$ will be the aggregate of all
the products of two unequal terms, the coefficient of $z^3$ will
be the aggregate of all products of three unequal terms,
and so on; from this it follows that
\[
R=1+\alpha z+\beta z^2+\gamma z^3+\delta z^4+\epsilon z^5+\textrm{etc.}
\]
according to the definitions given above (\S 3).

While if we put
\[
S=(1-az)(1-bz)(1-cz)(1-dz)(1-ez)\, \textrm{etc.},
\]
it will be just by making $z$ negative
\[
S=1-\alpha z+\beta z^2-\gamma z^3+\delta z^4-\epsilon z^5+\textrm{etc.}
\]

8. In order to compare these series $R$ and $S$ with the
preceding $P$ and $S$,
it should be noted that
\[
lR=l(1+az)+l(1+bz)+l(1+cz)+l(1+dz)+\textrm{etc.},
\]
which by taking the differentials will be
\[
\frac{dR}{Rdz}=\frac{a}{1+az}+\frac{b}{1+bz}+\frac{c}{1+cz}+\frac{d}{1+dz}+\textrm{etc.},
\]
which multiplied by $z$ gives the same previous expression which above we
called
$Q$, so that it will be
\[
Q=\frac{zdR}{Rdz}.
\]

Also in a similar way it will be
\[
\frac{dS}{Sdz}=-\frac{a}{1-az}-\frac{b}{1-bz}-\frac{c}{1-cz}-\textrm{etc.},
\]
from which one obtains
\[
P=\frac{-zdS}{Sdz}.
\]

9. Now, since 
\[
R=1+\alpha z+\beta z^2+\gamma z^3+\textrm{etc.},
\]
it will be
\[
\frac{zdR}{dz}=\alpha z+2\beta z^2+3\gamma z^3+4\delta z^4+5\epsilon z^5+\textrm{etc.}
\]
and hence
\begin{eqnarray*}
Q&=&Az-Bz^2+Cz^3-Dz^4+Ez^5-\textrm{etc.}\\
&=&\frac{\alpha z+2\beta z^2+3\gamma z^3+4\delta z^4+5\epsilon z^5+\textrm{etc.}}{1+\alpha z+\beta z^2+\gamma z^3+\delta z^4+\epsilon z^5+\textrm{etc.}}.
\end{eqnarray*}
Then, from the equality of these expressions
 the following relations between the letters
$A,B,C,D$ etc. and $\alpha,\beta,\gamma,\delta,\epsilon$ etc. occur:
\begin{eqnarray*}
A&=&\alpha,\\
\alpha A-B&=&2\beta,\\
\beta A-\alpha B+C&=&3\gamma,\\
\gamma A-\beta B-\alpha C-D&=&4\delta,\\
\delta A-\gamma B+\beta C-\alpha D+E&=&5\epsilon\\
\textrm{etc.}&&
\end{eqnarray*}

Indeed in a similar way, from the other equation $P=\frac{-zdS}{Sdz}$ it follows that
\begin{eqnarray*}
P&=&Az+Bz^2+Cz^3+Dz^4+Ez^5+\textrm{etc.}\\
&=&\frac{\alpha z-2\beta z^2+3\gamma z^3-4\delta z^4+5\delta z^5-\textrm{etc.}}{1-\alpha z+\beta z^2-\gamma z^3+\delta z^4-\epsilon z^5+\textrm{etc.}},
\end{eqnarray*}
which even yields the same determinations which we gave
above (\S 5). 

10. Also, by integrating the equation
$Q=\frac{zdR}{Rdz}$ it follows that $\int \frac{Qdz}{z}=lR$.
Indeed because $Q=Az-Bz^2+Cz^3-Dz^4+$ etc., it will be
\[
\int \frac{Qdz}{z}=Az-\frac{Bz^2}{2}+\frac{Cz^3}{3}-\frac{Dz^4}{4}+\textrm{etc.},
\]
and the value of this series thus expresses the logarithm of this series
\[
R=1+\alpha z+\beta z^2+\gamma z^3+\delta z^4+\textrm{etc.}
\]

Therefore because
\[
l(1+\alpha z+\beta z^2+\gamma z^3+\textrm{etc.})=Az-\frac{1}{2}Bz^2+\frac{1}{3}Cz^3-\frac{1}{4}Dz^4+\textrm{etc.},
\]
then from the equation $\int \frac{Pdz}{z}=-lS$ it will be
\[
l(1-\alpha z+\beta z^2-\gamma z^3+\textrm{etc.})=-Az-\frac{1}{2}Bz^2-\frac{1}{3}Cz^3-\frac{1}{4}Dz^4-\textrm{etc.}
\]
So if $k$ is written as the number whose hyperbolic logarithm is $=1$,
we will have
\[
1+\alpha z+\beta z^2+\gamma z^3+\delta z^4+\textrm{etc.}=k^{Az-\frac{1}{2}Bz^2+\frac{1}{3}Cz^3-\frac{1}{4}Dz^4+\textrm{etc.}}
\]
and
\[
1-\alpha z+\beta z^2-\gamma z^3+\delta z^4-\textrm{etc.}=k^{-Az-\frac{1}{2}Bz^2-\frac{1}{3}Cz^3-\frac{1}{4}Dz^4-\textrm{etc.}}.
\]

11. Also noteworthy are the reciprocals expressions of $R$ and $S$,
of course $\frac{1}{R}$ and $\frac{1}{S}$. Indeed it is
\[
\frac{1}{S}=\frac{1}{(1-az)(1-bz)(1-cz)(1-dz)\, \textrm{etc.}};
\]
to express the value of this fraction by a series whose
terms proceed according to powers of $z$, it is clear that all
these geometric progressions should be multiplied into each other
\begin{eqnarray*}
\frac{1}{1-az}&=&1+az+a^2z^2+a^3z^3+a^4z^4+\textrm{etc.},\\
\frac{1}{1-bz}&=&1+bz+b^2z^2+b^3z^3+b^4z^4+\textrm{etc.},\\
\frac{1}{1-cz}&=&1+cz+c^2z^2+c^3z^3+c^4z^4+\textrm{etc.},\\
\frac{1}{1-dz}&=&1+dz+d^2z^2+d^3z^3+d^4z^4+\textrm{etc.}\\
&&\textrm{etc.}
\end{eqnarray*}
In the product, after the first term $1$, the coefficient of $z$ will be
the sum of the quantities $a+b+c+d+\textrm{etc.}$, the coefficient
of $z^2$ will be the sum of the factors from two, not
excluding equal factors in the same product,
the coefficient of $z^3$ will be the sum of the factors from two,
and so on. We designated these sums of products above (\S 4)
with the letters of the Germanic alphabet
$\mathfrak{A},\mathfrak{B},\mathfrak{C},\mathfrak{D},\mathfrak{E}$ etc.
With these letters thus introduced, we will have
\[
\frac{1}{S}=1+\mathfrak{A}z+\mathfrak{B}z^2+\mathfrak{C}z^3+\mathfrak{D}z^4+\mathfrak{E}z^5+\textrm{etc.}
\]
and by treating the value of $R$ in a similar way it will be
\[
\frac{1}{R}=1-\mathfrak{A}z+\mathfrak{B}z^2-\mathfrak{C}z^3+\mathfrak{D}z^4-\mathfrak{E}z^5+\textrm{etc.}
\]

12. So these series are the reciprocals of those which we defined
under the letters $R$ and $S$ above (\S 7). And because of this it will be
\[
1=(1+\alpha z+\beta z^2+\gamma z^3+\delta z^4+\textrm{etc.})(1-\mathfrak{A}z+\mathfrak{B}z^2-\mathfrak{C}z^3+\mathfrak{D}z^4+\textrm{etc.})
\]
and even
\[
1=(1-\alpha z+\beta z^2-\gamma z^3+\delta z^4-\textrm{etc.})(1+\mathfrak{A}z+\mathfrak{B}z^2+\mathfrak{C}z^3+\mathfrak{D}z^4+\textrm{etc.}).
\]
From either one of these follows the same relation between the values
of the letters $\mathfrak{A},\mathfrak{B},\mathfrak{C},\mathfrak{D}$ etc.
and $\alpha,\beta,\gamma,\delta$ etc.; namely it will be
\begin{eqnarray*}
\mathfrak{A}-\alpha&=&0,\\
\mathfrak{B}-\alpha \mathfrak{A}+\beta&=&0,\\
\mathfrak{C}-\alpha \mathfrak{B}+\beta \mathfrak{A}-\gamma&=&0,\\
\mathfrak{D}-\alpha \mathfrak{C}+\beta \mathfrak{B}-\gamma \mathfrak{A}+\delta
&=&0\\
\textrm{etc.},&&
\end{eqnarray*}
which is the same relation we already gave above (\S 5).

13. But if we put $\frac{1}{R}=T$ and $\frac{1}{S}=V$, so that
\[
T=1-\mathfrak{A}z+\mathfrak{B}z^2-\mathfrak{C}z^3+\mathfrak{D}z^4-\textrm{etc.}
\]
and
\[
V=1+\mathfrak{A}z+\mathfrak{B}z^2+\mathfrak{C}z^3+\mathfrak{D}z^4+\textrm{etc.}
\]
it will be
\[
\frac{dR}{R}=-\frac{dT}{T} \quad \textrm{and} \quad \frac{dS}{S}=-\frac{dV}{V}
\]
and therefore it will become
\[
P=\frac{zdV}{Vdz} \quad \textrm{and} \quad Q=-\frac{zdT}{Tdz}.
\]
Now since 
\[
\frac{zdV}{dz}=\mathfrak{A}z+2\mathfrak{B}z^2+3\mathfrak{C}z^3+4\mathfrak{D}z^4+\textrm{etc.}
\]
and
\[
-\frac{zdT}{dz}=\mathfrak{A}z-2\mathfrak{B}z^2+3\mathfrak{C}z^3-4\mathfrak{D}z^4+\textrm{etc.},
\]
by writing the appropriate values from \S 6 in place of $P$ and $Q$ we will
have these equations
\[
Az+Bz^2+Cz^3+Dz^4+\textrm{etc.}=\frac{\mathfrak{A}z+2\mathfrak{B}z^2+3\mathfrak{C}z^3+4\mathfrak{D}z^4+\textrm{etc.}}{1+\mathfrak{A}z+\mathfrak{B}z^2+\mathfrak{C}z^3+\mathfrak{D}z^4+\textrm{etc.}}
\]
and
\[
Az-Bz^2+Cz^3-Dz^4+\textrm{etc.}=\frac{\mathfrak{A}z-2\mathfrak{B}z^2+3\mathfrak{C}z^3
-4\mathfrak{D}z^4+\textrm{etc.}}{1-\mathfrak{A}z
+\mathfrak{B}z^2-\mathfrak{C}z^3+\mathfrak{D}z^4-\textrm{etc.}},
\]
from which the same relation between the letters $A,B,C,D$ etc. and
$\mathfrak{A},\mathfrak{B},\mathfrak{C},\mathfrak{D}$ etc. follows
as we gave above (\S 5). Namely it will be
\begin{eqnarray*}
\mathfrak{A}&=&A,\\
2\mathfrak{B}&=&\mathfrak{A}A+B,\\
3\mathfrak{C}&=&\mathfrak{B}A+\mathfrak{A}B+C,\\
4\mathfrak{D}&=&\mathfrak{C}A+\mathfrak{B}B+\mathfrak{A}C+D,\\
5\mathfrak{E}&=&\mathfrak{D}A+\mathfrak{C}B+\mathfrak{B}C+\mathfrak{A}D+E\\
&&\textrm{etc.}
\end{eqnarray*}

14. From the equations given in \S 12 it follows that
\[
l(1+\alpha z+\beta z^2+\gamma z^3+\textrm{etc.})=-l(1-\mathfrak{A}z+\mathfrak{B}z^2-\mathfrak{C}z^3+\textrm{etc.})
\]
and
\[
l(1-\alpha z+\beta z^2-\gamma z^3+\textrm{etc.})=-l(1+\mathfrak{A}z+\mathfrak{B}z^2+\mathfrak{C}z^3+\textrm{etc.}).
\]
Then by applying these to \S 10 it will be
\[
l(1-\mathfrak{A}z+\mathfrak{B}z^2-\mathfrak{C}z^3+\textrm{etc.})=-Az+\frac{1}{2}Bz^2
-\frac{1}{3}Cz^3+\frac{1}{4}Dz^4-\textrm{etc.}
\]
and
\[
l(1+\mathfrak{A}z+\mathfrak{B}z^2+\mathfrak{C}z^3+\textrm{etc.})=Az+
\frac{1}{2}Bz^2+\frac{1}{3}Cz^3+\frac{1}{4}Dz^4+\textrm{etc.}
\]
And hence by taking $k$ as the number whose logarithm is $=1$, it will be
\[
1-\mathfrak{A}z+\mathfrak{B}z^2-\mathfrak{C}z^3+\textrm{etc.}=k^{-Az+\frac{1}{2}Bz^2-\frac{1}{3}Cz^3+\frac{1}{4}Dz^4-\textrm{etc.}}
\]
and
\[
1+\mathfrak{A}z+\mathfrak{B}z^2+\mathfrak{C}z^3+\textrm{etc.}=k^{Az+\frac{1}{2}Bz^2+\frac{1}{3}Cz^3+\frac{1}{4}Dz^4+\textrm{etc.}}.
\]

15. Now if the letters $R$ and $S$ retain the values assumed above (\S 7), it will be 
\begin{eqnarray*}
1+\alpha z+\beta z^2+\gamma z^3+\delta z^4+\textrm{etc.}&=&R,\\
1-\mathfrak{A}z+\mathfrak{B}z^2-\mathfrak{C}z^3+\mathfrak{D}z^4-\textrm{etc.}&=&\frac{1}{R}
\end{eqnarray*}
and
\begin{eqnarray*}
1-\alpha z+\beta z^2-\gamma z^3+\delta z^4-\textrm{etc.}&=&S,\\
1+\mathfrak{A}z+\mathfrak{B}z^2+\mathfrak{C}z^3+\mathfrak{D}z^4+\textrm{etc.}&=&\frac{1}{S}.
\end{eqnarray*}
From these the following conclusions are deduced
\begin{eqnarray*}
1+\beta z^2+\delta z^4+\zeta z^6+\theta z^8+\textrm{etc.}&=&\frac{R+S}{2},\\
\alpha z+\gamma z^3+\epsilon z^5+\eta z^7+\iota z^9+\textrm{etc.}&=&\frac{R-S}{2},\\
1+\mathfrak{B}z^2+\mathfrak{D}z^4+\mathfrak{F}z^6+\mathfrak{H}z^8+\textrm{etc.}&=&\frac{R+S}{2RS},\\
\mathfrak{A}z+\mathfrak{C}z^3+\mathfrak{E}z^5+\mathfrak{G}z^7+\mathfrak{I}z^9+\textrm{etc.}&=&\frac{R-S}{2RS}
\end{eqnarray*}
and hence this proportion is obtained
\[
\begin{split}
&1+\beta z^2+\delta z^4+\zeta z^6+\textrm{etc.}:\alpha z+\gamma z^3+\epsilon z^5+\eta z^7+\textrm{etc.}\\
=&1+\mathfrak{B}z^2+\mathfrak{D}z^4+\mathfrak{F}z^6+\textrm{etc.}:\mathfrak{A}z+\mathfrak{C}z^3+\mathfrak{E}z^5+\mathfrak{G}z^7+\textrm{etc.}
\end{split}
\]
Since it is also
\begin{eqnarray*}
R-1&=&\alpha z+\beta z^2+\gamma z^3+\delta z^4+\textrm{etc.},\\
1-\frac{1}{R}&=&\mathfrak{A}z-\mathfrak{B}z^2+\mathfrak{C}z^3-\mathfrak{D}z^4+\textrm{etc.},
\end{eqnarray*}
it will be
\[
R=\frac{\alpha z+\beta z^2+\gamma z^3+\delta z^4+\textrm{etc.}}{\mathfrak{A}z-\mathfrak{B}z^2+\mathfrak{C}z^3-\mathfrak{D}z^4+\textrm{etc.}},
\]
and in a similar way, because
\begin{eqnarray*}
1-S&=&\alpha z-\beta z^2+\gamma z^3-\delta z^4+\textrm{etc.},\\
\frac{1}{S}-1&=&\mathfrak{A}z+\mathfrak{B}z^2+\mathfrak{C}z^3+\mathfrak{D}z^4+\textrm{etc.}
\end{eqnarray*}
it will be
\[
S=\frac{\alpha z-\beta z^2+\gamma z^3-\delta z^4+\textrm{etc.}}{\mathfrak{A}z+\mathfrak{B}z^2+\mathfrak{C}z^3+\mathfrak{D}z^4+\textrm{etc.}}.
\]

16. Next indeed if as above (\S 6) we put
\begin{eqnarray*}
P&=&Az+Bz^2+Cz^3+Dz^4+\textrm{etc.},\\
Q&=&Az-Bz^2+Cz^3-Dz^4+\textrm{etc.},
\end{eqnarray*}
it will be from paragraph 9
\begin{eqnarray*}
\alpha z+2\beta z^2+3\gamma z^3+4\delta z^4+\textrm{etc.}&=&QR,\\
\alpha z-2\beta z^2+3\gamma z^3-4\delta z^4+\textrm{etc.}&=&PS
\end{eqnarray*}
and in a similar way from paragraph 13 we will have\footnote{Translator: The {\em Opera omnia},
Series I, Volume 2, p. 168 and the original version in the {\em Commentarii},
Volume 13, p. 74 both have $+4\mathfrak{D}z^4$ on the second line, which
is
incorrect.}
\begin{eqnarray*}
\mathfrak{A}z+2\mathfrak{B}z^2+3\mathfrak{C}z^3+4\mathfrak{D}z^4+\textrm{etc.}&=&\frac{P}{S},\\
\mathfrak{A}z-2\mathfrak{B}z^2+3\mathfrak{C}z^3-4\mathfrak{D}z^4+
\textrm{etc.}&=&\frac{Q}{R}.
\end{eqnarray*}
From these the following corollaries are easily derived:
\[
\frac{\alpha z-2\beta z^2+3\gamma z^3-4\delta z^4+\textrm{etc.}}{Az+Bz^2+Cz^3+Dz^4+\textrm{etc.}}=S=\frac{Az+Bz^2+Cz^3+Dz^4+\textrm{etc.}}{\mathfrak{A}z+2\mathfrak{B}z^2+3\mathfrak{C}z^3+4\mathfrak{D}z^4+\textrm{etc.}},
\]
\[
\frac{\alpha z+2\beta z^2+3\gamma z^3+4\delta z^4+\textrm{etc.}}{Az-Bz^2+Cz^3-Dz^4+\textrm{etc.}}=R=\frac{Az-Bz^2+Cz^3-Dz^4+\textrm{etc.}}{\mathfrak{A}z-2\mathfrak{B}z^2+3\mathfrak{C}z^3-4\mathfrak{D}z^4+\textrm{etc.}},
\]
For the letters $R$ and $S$ we have these quintuple values
\begin{eqnarray*}
R&=&1+\alpha z+\beta z^2+\gamma z^3+\delta z^4+\textrm{etc.},\\
R&=&\frac{1}{1-\mathfrak{A}z+\mathfrak{B}z^2-\mathfrak{C}z^3+\mathfrak{D}z^4-\textrm{etc.}},\\
R&=&\frac{\alpha z+\beta z^2+\gamma z^3+\delta z^4+\textrm{etc.}}{\mathfrak{A}z-\mathfrak{B}z^2+\mathfrak{C}z^3-\mathfrak{D}z^4+\textrm{etc.}},\\
R&=&\frac{\alpha z+2\beta z^2+3\gamma z^3+4\delta z^4+\textrm{etc.}}{Az-Bz^2+Cz^3-4z^4+\textrm{etc.}},\\
R&=&\frac{Az-Bz^2+Cz^3-Dz^4+\textrm{etc.}}{\mathfrak{A}z-2\mathfrak{B}z^2+3\mathfrak{C}z^3-4\mathfrak{D}z^4+\textrm{etc.}},
\end{eqnarray*}
in which putting $-z$ in place of $z$ everywhere yields the values for
$S$. 
And from varied combinations of these
five values, a great number of properties can be elicited 
which the three orders of our letters $A,B,C,D$ etc.,
$\alpha,\beta,\gamma,\delta$ etc., $\mathfrak{A},\mathfrak{B},
\mathfrak{C},\mathfrak{D}$ etc.
hold between each other,
which however we will refrain from pursuing now.

17. From this it is quite clear that that we can descend from what has
been explained already to more particular cases, and first indeed this
infinite geometric progression will be taken for the series of letters
$a,b,c,d$ etc.
\[
n,n^2,n^3,n^4,n^5,n^6 \, \textrm{etc.};
\]
with these successively introduced into the above formulae we will have:
\begin{eqnarray*}
A&=&n+n^2+n^3+n^4+n^5+\textrm{etc.}=\frac{n}{1-n},\\
B&=&n^2+n^4+n^6+n^8+n^{10}+\textrm{etc.}=\frac{nn}{1-nn},\\
C&=&n^3+n^6+n^9+n^{12}+n^{15}+\textrm{etc.}=\frac{n^3}{1-n^3},\\
D&=&n^4+n^8+n^{12}+n^{16}+n^{20}+\textrm{etc.}=\frac{n^4}{1-n^4}\\
&&\textrm{etc.}
\end{eqnarray*}
Now from \S 6, we will get two values for the letters $P$ and $Q$, which
will be
\begin{eqnarray*}
P&=&\frac{nz}{1-nz}+\frac{n^2z}{1-n^2z}+\frac{n^3z}{1-n^3z}+\frac{n^4z}{1-n^4z}+\textrm{etc.},\\
Q&=&\frac{nz}{1+nz}+\frac{n^2z}{1+n^2z}+\frac{n^3z}{1+n^3z}+\frac{n^4z}{1+n^4z}+\textrm{etc.}
\end{eqnarray*}
and then from the values found for the letters $A,B,C,D$ etc., these other
values will arise
\begin{eqnarray*}
P&=&\frac{nz}{1-n}+\frac{n^2z^2}{1-nn}+\frac{n^3z^3}{1-n^3}+\frac{n^4z^4}{1-n^4}+\textrm{etc.},\\
Q&=&\frac{nz}{1-n}-\frac{n^2z^2}{1-nn}+\frac{n^3z^3}{1-n^3}-\frac{n^4z^4}{1-n^4}+\textrm{etc.}
\end{eqnarray*}

18. Next, from paragraph 7 we will have the following expressions
for $R$ and $S$
\begin{eqnarray*}
R&=&(1+nz)(1+n^2z)(1+n^3z)(1+n^4z)\, \textrm{etc.},\\
S&=&(1-nz)(1-n^2z)(1-n^3z)(1-n^4z)\, \textrm{etc.}
\end{eqnarray*}
These factors actually multiplied into each other
and arranged according to dimensions of $z$ yield these series
for $R$ and $S$
\begin{eqnarray*}
R&=&1+\alpha z+\beta z^2+\gamma z^3+\delta z^4+\textrm{etc.},\\
S&=&1-\alpha z+\beta z^2-\gamma z^3+\delta z^4-\textrm{etc.},
\end{eqnarray*}
where the letters $\alpha,\beta,\gamma,\delta$ etc. are thus determined
from the supposed series $n,n^2,n^3,n^4,n^5,n^6,n^7$ etc.,
so that it will be:

I. $\alpha=$ the sum of all the terms; whence it will be 
\[
\alpha=n+n^2+n^3+n^4+n^5+n^6+n^7+\textrm{etc.},
\]
which is supposed to be a geometric progression,
in which each power of $n$ occurs and has the coefficient $+1$.

II. $\beta=$ the sum of the factors from two terms; whence it will be
\[
\beta=n^3+n^4+2n^5+2n^6+3n^7+3n^8+4n^9+4n^{10}+\textrm{etc.},
\]
in which series after the third power all the following powers
of $n$ occur;
moreover, each power occurs as often as it can be made
by multiplying two terms of the series $\alpha$.
Since however the multiplication of powers consists in the addition
of exponents, the coefficient of any power of $n$ will appear
in the series $\beta$ in as many ways as the exponent of $n$
can be distributed into two unequal parts, or as many ways
as this exponent $n$ can be produced from the addition
of two unequal integral numbers.
Thus the coefficient of the tenth power $n^{10}$ is $4$,
because $10$ can be distributed into two unequal parts in four
ways, namely,
\[
\begin{array}{ll}
10=1+9,&10=3+7,\\
10=2+8,&10=4+6.
\end{array}
\]

III. $\gamma=$ the sum of the factors from three unequal terms of
the series $\alpha$;
whence it will be
\[
\gamma=n^6+n^7+2n^8+3n^9+4n^{10}+5n^{11}+7n^{12}+8n^{13}+\textrm{etc.},
\]
in which after the sixth power, all the following
powers of $n$ occur. Moreover the coefficient of each power indicates
how many ways the exponent can be distributed into three unequal parts,
or as often as the same exponent can be produced from the addition
of three mutually unequal integral numbers.
Thus the power $n^{12}$ has the coefficient $7$, because the exponent $12$
can be partitioned into three unequal parts in seven ways, as
\[
\begin{array}{ll}
12=1+2+9,&12=1+5+6,\\
12=1+3+8,&12=2+3+7,\\
12=1+4+7,&12=2+4+6,
\end{array}
\]
\[
12=3+4+5.
\]

IV. $\delta=$ the sum of the factors from four mutually unequal
terms of the series $\alpha$; whence it will be
\[
\delta=n^{10}+n^{11}+2n^{12}+3n^{13}+5n^{14}+6n^{15}+9n^{16}+\textrm{etc.},
\]
whose first power is $n^{10}$, 
whose exponent is clearly $1+2+3+4$, or the fourth trigonal number.
Each of the following powers appears as often as its exponent 
can be made from the addition of four mutually unequal integral numbers.
Thus the sixteenth power $n^{16}$ has the coefficient $9$,
because $16$ 
can be distributed into four mutually unequal parts in nine ways. 
These nine partitions are
\[
\begin{array}{ll}
16=1+2+3+10,&16=1+3+4+8,\\
16=1+2+4+9,&16=1+3+5+7,\\
16=1+2+5+8,&16=1+4+5+6,\\
16=1+2+6+7,&16=2+3+4+7,
\end{array}
\]
\[
16=2+3+5+6.
\]
The same kind of thing happens for the values of the following letters
$\epsilon,\zeta,\eta$ etc., 
which will be
\begin{eqnarray*}
\epsilon&=&n^{15}+n^{16}+2n^{17}+3n^{18}+5n^{19}+7n^{20}+10n^{21}+\textrm{etc.},\\
\zeta&=&n^{21}+n^{22}+2n^{23}+3n^{24}+5n^{25}+7n^{26}+11n^{27}+\textrm{etc.},\\
\eta&=&n^{28}+n^{29}+2n^{30}+3n^{31}+5n^{32}+7n^{33}+11n^{34}+\textrm{etc.},\\
&&\textrm{etc.}
\end{eqnarray*}
In all these series, the coefficient of each power of $n$ indicates
how many different ways the exponent of $n$ can be resolved into as
many 
unequal parts as the series is numbered from the first.
In other words, the coefficient of any term indicates how many ways the
exponent of $n$ can be made from the addition of as many mutually
unequal integral numbers as the position of the series from which the
term is taken is numbered starting at $\alpha$.
Thus in the seventh series the coefficient of the power $n^{34}$ is $11$,
because the number $34$ can be distributed into seven unequal parts in seven ways; these distributions are
\begin{eqnarray*}
34&=&1+2+3+4+5+6+13,\\
34&=&1+2+3+4+5+7+12,\\
34&=&1+2+3+4+5+8+11,\\
34&=&1+2+3+4+5+9+10,\\
34&=&1+2+3+4+6+7+11,\\
34&=&1+2+3+4+6+8+10,\\
34&=&1+2+3+4+7+8+9,\\
34&=&1+2+3+5+6+7+10,\\
34&=&1+2+3+5+6+8+9,\\
34&=&1+2+4+5+6+7+9,\\
34&=&1+3+4+5+6+7+8.
\end{eqnarray*}
And indeed from these, the nature of the series
which appear for the letters
$\alpha,\beta,\gamma,\delta$ etc.
is easily seen.

19. Therefore for investigating how many different ways a number
can be distributed into a given number of unequal parts, the series expressing the letters
$\alpha,\beta,\gamma,\delta$ etc. can be formed, although this work
will be exceedingly tiresome. In turn, however, assuming these series as
known and already formed, a not inelegant problem can be solved, which was
thus proposed to me by the Insightful Naud\'e:

{\em To define how many ways a given number  can be produced from the addition
of several integral numbers, mutually unequal, the number of which is given.}

The most Insightful Proposer searched thus for how many different ways the number
$50$ can arise from the addition of seven unequal integral numbers.
For resolving this question it is clear that the appropriate series to consider
is $\eta$,
in which the coefficient of any term indicates how many different ways the exponent of $n$ can be resolved into $7$ unequal parts. Hence the series
\[
\eta=n^{28}+n^{29}+2n^{30}+3n^{31}+5n^{32}+7n^{33}+11n^{34}+\textrm{etc.}
\]
should be continued to the term at which the
fiftieth power of $n$ is contained,
whose coefficient, which will be $522$,
shows that the number $50$ can be produced in altogether $522$ different
ways from the addition of seven mutually unequal integral numbers.
Thus it is clear that if a convenient and simple way were obtained of forming
these series $\alpha,\beta,\gamma,\delta$ etc., Naud\'e's problem would
be brought to a most perfect solution.

20. Thus since a way was given above (\S 5 and 9) 
for finding the values of the letters
$\alpha,\beta,\gamma,\delta$ etc.
from the known values of the letters $A,B,C,D$ etc.,
for the present question we can easily obtain a solution,
since from \S 17 we have the known values 
$A,B,C,D$ etc.; and thus it follows that
\begin{eqnarray*}
\alpha&=&A,\\
\beta&=&\frac{\alpha A-B}{2},\\
\gamma&=&\frac{\beta A-\alpha B+C}{3},\\
\delta&=&\frac{\gamma A-\beta B+\alpha C-D}{4},\\
\epsilon&=&\frac{\delta A-\gamma B+\beta C-\alpha D+E}{5}\\
&&\textrm{etc.}
\end{eqnarray*}
We therefore obtain from these
\begin{eqnarray*}
\alpha&=&\frac{n}{1-n},\\
2\beta&=&\frac{\alpha n}{1-n}-\frac{nn}{1-nn},\\
3\gamma&=&\frac{\beta n}{1-n}-\frac{\alpha n^2}{1-n^2}+\frac{n^3}{1-n^3},\\
4\delta&=&\frac{\gamma n}{1-n}-\frac{\beta n^2}{1-n^2}+\frac{\alpha n^3}{1-n^3}
-\frac{n^4}{1-n^4}\\
&&\textrm{etc.}
\end{eqnarray*}
Now if however in place of 
$\alpha,\beta,\gamma$ etc. the
previous values are successively 
substituted, it follows
\begin{eqnarray*}
\alpha&=&\frac{n}{1-n},\\
\beta&=&\frac{n^3}{(1-n)(1-nn)},\\
\gamma&=&\frac{n^6}{(1-n)(1-nn)(1-n^3)},\\
\delta&=&\frac{n^{10}}{(1-n)(1-n^2)(1-n^3)(1-n^4)},\\
\epsilon&=&\frac{n^{15}}{(1-n)(1-n^2)(1-n^3)(1-n^4)(1-n^5)}\\
&&\textrm{etc.}
\end{eqnarray*}
Thus we can understand that, in this case,
\begin{eqnarray*}
\alpha&=&A,\\
\beta&=&AB,\\
\gamma&=&ABC,\\
\delta&=&ABCD,\\
\epsilon&=&ABCDE\\
&&\textrm{etc.}
\end{eqnarray*}

21. The law according to which the values of the letters $\alpha,\beta,\gamma,\delta$
etc. have been found, whose truth is observed by expanding several
formulae, is so far not apparent except by induction. So its truth
may thus be securely confirmed, it will be convenient to elicit
the same law of the progression in a totally different way, in which induction
plays no part. Now, let it be proposed to investigate the values
of the letters $\alpha,\beta,\gamma,\delta$ etc. which occur
in the series
\[
R=1+\alpha z+\beta z^2+\gamma z^3+\delta z^4+\epsilon z^5+\textrm{etc.}
\]
If as we had assumed initially
\[
R=(1+nz)(1+n^2z)(1+n^3z)(1+n^4z)\cdots,
\]
it should be noted that if $nz$ is written in place of $z$, the expression
which just now $R$ was equal to is changed into the form
\[
(1+n^2z)(1+n^3z)(1+n^4z)(1+n^5z)\cdots,
\]
which multiplied by $1+nz$ produces the prior expression. Therefore
we can rightly conclude that if in the series
\[
1+\alpha z+\beta z^2+\gamma z^3+\delta z^4+\epsilon z^5+\textrm{etc.}
\]
we write $nz$ in place of $z$, so that we have
\[
1+\alpha nz+\beta n^2z^2+\gamma n^3z^3+\delta n^4z^4+\epsilon n^5z^5+\textrm{etc.},
\]
and we then multiply this expression by $1+nz$, then the product, which will
be
\[
\begin{array}{rrrrrrr}
1&+\alpha nz&+\beta n^2z^2&+\gamma n^3z^3&+\delta n^4z^4&+\epsilon n^5z^5&+\textrm{etc.}\\
&+nz&+\alpha n^2z^2&+\beta n^3z^3&+\gamma n^4z^4&+\delta n^5z^5&+\textrm{etc.},
\end{array}
\]
should be equal to that of the prior series
\[
1+\alpha z+\beta z^2+\gamma z^3+\delta z^4+\epsilon z^5+\textrm{etc.}
\]
But if we equate the coefficients of the corresponding terms, we will
obtain the following determinations for $\alpha,\beta,\gamma$ etc.
\[
\begin{split}
&\alpha=\frac{n}{1-n}=\frac{n}{1-n},\\
&\beta=\frac{\alpha n^2}{1-n^2}=\frac{n^3}{(1-n)(1-n^2)},\\
&\gamma=\frac{\beta n^3}{1-n^3}=\frac{n^6}{(1-n)(1-n^2)(1-n^3)},\\
&\delta=\frac{\gamma n^4}{1-n^4}=\frac{n^{10}}{(1-n)(1-n^2)(1-n^3)(1-n^4)}\\
&\textrm{etc.}
\end{split}
\]

22. In this way we have thus found convenient enough expressions for
the sums of these series $\alpha,\beta,\gamma,\delta$ etc.,
from which in turn these series can be formed.
Namely since this series proceed according to powers of $n$,
 this will arise if the expressions of these sums are expanded by division
in the usual way into infinite series proceeding according to powers of $n$.
Let this division be done, 
and it is
clear that all the series $\alpha,\beta,\gamma,\delta$ etc. belong to
a type which is commonly referred to by the name of recurrent series; they
are such that any term can be determined from several of the preceding.
So it will be clear how in each of these series any term is formed
from the preceding, let us expand the denominators of
these expressions found for the letters $\alpha,\beta,\gamma,\delta$ etc.
by actual multiplication, and having done this we will have
\begin{eqnarray*}
\alpha&=&\frac{n}{1-n},\\
\beta&=&\frac{n^3}{1-n-n^2+n^3},\\
\gamma&=&\frac{n^6}{1-n-n^2+n^4+n^5-n^6},\\
\delta&=&\frac{n^{10}}{1-n-n^2+2n^5-n^8-n^9+n^{10}},\\
\epsilon&=&\frac{n^{15}}{1-n-n^2+n^5+n^6+n^7-n^8-n^9-n^{10}+n^{13}+n^{14}
-n^{15}},\\
\zeta&=&\frac{n^{21}}{1-n-n^2+n^5+2n^7-n^9-n^{10}-n^{11}-n^{12}
+2n^{14}+n^{16}-n^{19}-n^{20}+n^{21}}\\
&&\textrm{etc.}
\end{eqnarray*}
And from these denominators we understand how in each series any term
is composed from the preceding, if the precept by which these recurrent
series are formed is called upon.

23. And from the form of the expressions found for the letters
$\alpha,\beta,\gamma,\delta$ etc., by which each is the product of the preceding
and some new factor, 
another method 
is deduced which is suitable enough for finding the next series from
any series which has already been found.
Thus, since the series $\alpha=\frac{n}{1-n}$ is a geometric progression
\[
\alpha=n+n^2+n^3+n^4+n^5+n^6+n^7+\textrm{etc.},
\]
from which comes the series $\beta$, if it is multiplied by $\frac{n^2}{1-n^2}$ or if 
it is multiplied by this geometric progression
\[
n^2+n^4+n^6+n^8+n^{10}+n^{12}+n^{14}+\textrm{etc.}
\]
Next, with the series $\beta$ having been solidly found, if it is multiplied by
\[
\frac{n^3}{1-n^3}=n^3+n^6+n^9+n^{12}+n^{15}+n^{18}+\textrm{etc.},
\]
series $\gamma$ is produced.
And this multiplied by
\[
\frac{n^4}{1-n^4}=n^4+n^8+n^{12}+n^{16}+n^{20}+n^{24}+\textrm{etc.}
\]
yields the series $\delta$.
And by so on multiplying each series of this order by a certain geometric progression the following series comes out. In this manner, these series can be continued
as far as we want; and thus the above problem proposed by the Insightful Naud\'e is resolved.

24. And each of the series will be able to be easily continued by means
of the preceding, if we consider the way in which the value of each
of the letters $\alpha,\beta,\gamma,\delta$ etc. is determined from the
preceding.
Thus, since $\beta=\frac{\alpha n^2}{1-n^2}$, it will be $\beta=\beta nn+\alpha nn$;
therefore if to the series $\beta$ multiplied by $nn$ is added
the series $\alpha$ multiplied by $nn$, the series $\beta$ will be created.
Then, since it is clear that the first term of the series $\beta$ is $n^3$,
let us put
\[
\beta=\mathfrak{a}n^3+\mathfrak{b}n^4+\mathfrak{c}n^5+\mathfrak{d}n^6
+\mathfrak{e}n^7+\mathfrak{f}n^8+\mathfrak{g}n^9+\textrm{etc.}
\]
and it will be
\[
\begin{array}{rcllllllll}
\beta n^2&=&&&+\mathfrak{a}n^5&+\mathfrak{b}n^6&+\mathfrak{c}n^7&+\mathfrak{d}n^8&+\mathfrak{e}n^9&+\textrm{etc.},\\
\alpha n^2&=&n^3&+n^4&+n^5&+n^6&+n^7&+n^8&+n^9&+\textrm{etc.}
\end{array}
\]
Now with the terms equated, because $\beta=\beta nn+\alpha nn$ we will have
\[
\begin{array}{ll}
\mathfrak{a}=1,&\mathfrak{e}=\mathfrak{c}+1=3,\\
\mathfrak{b}=1,&\mathfrak{f}=\mathfrak{d}+1=3,\\
\mathfrak{c}=\mathfrak{a}+1=2,&\mathfrak{g}=\mathfrak{e}+1=4,\\
\mathfrak{d}=\mathfrak{b}+1=2,&\mathfrak{h}=\mathfrak{f}+1=4
\end{array}
\]
\[
\textrm{etc.}
\]
In a similar way, since $\gamma=\frac{\beta n^3}{1-n^3}$ or
$\gamma=\gamma n^3+\beta n^3$, the series $\gamma$ will be formed from the
series $\beta$,
and in turn the series $\delta$ is produced from the series $\gamma$, by means of the
equation $\delta=\delta n^4+\gamma n^4$; and all the following
will be  dispatched likewise.

25. Because in the expression
\[
R=1+\alpha z+\beta z^2+\gamma z^3+\delta z^4+\textrm{etc.}
\]
we have found the values of the letters $\alpha,\beta,\gamma,\delta$ etc.
and it is
\[
R=(1+nz)(1+n^2z)(1+n^3z)(1+n^4z)\ldots,
\]
this product, apparently
\[
(1+nz)(1+n^2z)(1+n^3z)(1+n^4z)\ldots,
\]
will be converted from the infinite factors into this series proceeding
according to powers of $z$
\[
1+\frac{nz}{1-n}+\frac{n^3z^3}{(1-n)(1-n^2)}
+\frac{n^6z^3}{(1-n)(1-n^2)(1-n^3)}
+\frac{n^{10}z^4}{(1-n)(1-n^2)(1-n^3)(1-n^4)}
+\textrm{etc.}
\]
And from \S 10 the hyperbolic logarithm of the sum of this series will be
\[
=\frac{nz}{1-n}-\frac{nnz^2}{2(1-n^2)}
+\frac{n^3z^3}{3(1-n^3)}
-\frac{n^4z^4}{4(1-n^4)}+\textrm{etc.}
\]
Or if $k$ is written for the number whose logarithm $=1$, it will be
\[
k^{\frac{nz}{1-n}-\frac{n^2z^2}{2(1-n^2)}+\frac{n^3z^3}{3(1-n^3)}
-\frac{n^4z^4}{4(1-n^4)}+\textrm{etc.}}=R
\]
or this exponential expression is equal to the sum of series which we made
changed into the value of $R$.

26. Truly to turn the proposed problem, which is to define how many
different ways a given number $m$ can be partitioned into $\mu$ integral
parts that are unequal to each other, let us indicate this number of ways
which we are seeking by the notation
\[
m^{(\mu)i},
\]
which from now on will indicate to us the number of ways in which the number
$m$ can be produced from the addition of $\mu$ mutually unequal integral
numbers; and for denoting the inequality of these parts we have adjoined
the letter $i$ above, which will be omitted if the question takes the form
of finding altogether the number of ways in which the given number $m$ can be
distributed into $\mu$ parts, either equal or unequal. 
After this, a solution to the problem can easily be shown

27. Thus this number of ways $m^{(\mu)i}$ will be the coefficient of the
power $n^m$ in the series $\alpha,\beta,\gamma,\delta,\epsilon$ etc., which
first at $\alpha$ are numbered as far as $\mu$ contains unities. The sum
of this series is
\[
=\frac{n^{\frac{\mu(\mu+1)}{1\cdot 2}}}{(1-n)(1-n^2)(1-n^3)(1-n^4)\cdots(1-n^\mu)}
\]
and hence the general term of the series which arises from this form is
$=m^{(\mu)i}n^m$. Moreover, the general term of the series which arises
from the form
\[
\frac{n^{\frac{\mu(\mu-1)}{1\cdot 2}}}{(1-n)(1-n^2)(1-n^3)(1-n^4)\cdots(1-n^\mu)}
\]
will be $=m^{(\mu)i}n^{m-\mu}$ or for the same power of $n$ the general
term will be $=(m+\mu)^{(\mu)i}n^m$.
The prior expression is subtracted from the latter, and the general
term of the remaining expression
\[\frac{n^{\frac{\mu(\mu-1)}{1\cdot 2}}}{(1-n)(1-n^2)(1-n^3)(1-n^4)\cdots(1-n^{\mu-1})}\]
will be $=n^m((m+\mu)^{(\mu)i}-m^{(\mu)i})$;
moreover the general term of this series is $m^{(\mu-1)i}n^m$, whence we will
have
\[m^{(\mu-1)i}=(m+\mu)^{(\mu)i}-m^{(\mu)i},\]
from which we arrive at the rule that
\[
(m+\mu)^{(\mu)i}=m^{(\mu)i}+m^{(\mu-1)i},
\]
by means of which, if the number of different ways the number $m$ can be distributed into
$\mu$ and $\mu-1$ unequal parts were known, by adding these two numbers
would follow the number of ways in which the larger number $m+\mu$ can be distributed into
$\mu$ unequal parts. And thus the resolution of more difficult cases
is reduced to simpler ones, and with these known finally to the simpler
ones; it is of course clear that if $m<\frac{\mu\mu+\mu}{2}$ then
$m^{(\mu)i}=0$, and if $m=\frac{\mu\mu+\mu}{2}$ then
it will be $m^{(\mu)i}=1$.

28. Since the formula $m^{(\mu)i}n^m$ is the general term of the
expression
\[
\frac{n^{\frac{\mu(\mu+1)}{2}}}{(1-n)(1-n^2)(1-n^3)\cdots(1-n^\mu)},
\]
let us see what kind of series this expression
\[
\frac{1}{(1-n)(1-n^2)(1-n^3)\cdots(1-n^\mu)},
\]
produces
if expanded and arranged according
to the dimensions of $n$. Let us put it to produce this series
\[
1+pn+qn^2+rn^3+sn^4+tn^5+\textrm{etc.},
\]
from whose generation it is clear that the coefficient of any power
of $n$ shows how many different ways the exponent of $n$ can be produced
by addition from the given numbers
\[
1, 2,3,4,5,6,\ldots, \mu;
\]
and here neither is a certain number of parts prescribed from which it is
formed, nor is the condition put that the must are unequal to each other.
Therefore the expression $m^{(\mu)i}$ will indicate altogether how many
ways the number $m-\frac{\mu(\mu+1)}{2}$ can be produced by addition
from the numbers $1,2,3,4,5,\ldots, \mu$.
Thus if one searches for how many different ways the number $50$
can be distributed into $7$ unequal parts, because $m=50$ and $\mu=7$
the question is thus reduced to investigating how many different
ways the number $50-28$ or $22$ can arise by addition from the seven numbers
$1,2,3,4,5,6,7$. 
With this understood, both varieties of this question can be resolved
in a single effort.

29. With the letters $\alpha,\beta,\gamma,\delta$ etc. defined for the
case where we have assumed the geometric progression $n,n^2,n^3,n^4,n^5$
etc. in place of the letters $a,b,c,d$ etc., 
order requires that we also inquire into the values of
the third order $\mathfrak{A},\mathfrak{B},\mathfrak{C},\mathfrak{D},
\mathfrak{E}$, etc.
But we have employed the letters $\mathfrak{A},\mathfrak{B},\mathfrak{C},
\mathfrak{D}$ etc., with equal values, in the series $\frac{1}{R}$
and $\frac{1}{S}$;
for we assumed above (\S 11) that
\[
\frac{1}{S}=1+\mathfrak{A}z+\mathfrak{B}z^2+\mathfrak{C}z^3+
\mathfrak{D}z^4+\mathfrak{E}z^5+\textrm{etc.}
\]
and
\[
\frac{1}{R}=1-\mathfrak{A}z+\mathfrak{B}z^2-\mathfrak{C}z^3+\mathfrak{D}z^4
-\mathfrak{E}z^5+\textrm{etc.}
\]
where the original values taken for $R$ and $S$ were
\[
R=(1+nz)(1+n^2z)(1+n^3z)(1+n^4z)\, \textrm{etc.},
\]
\[
S=(1-nz)(1-n^2z)(1-n^3z)(1-n^4z)\, \textrm{etc.}
\]
It is apparent hence that the series $\frac{1}{S}=1+\mathfrak{A}z+
\mathfrak{B}z^2+\mathfrak{C}z^3+\mathfrak{D}z^4+$ etc. will arise if
these innumerable geometric
progressions are multiplied by each other
\begin{eqnarray*}
\frac{1}{1-nz}&=&1+nz+n^2z^2+n^3z^3+n^4z^4+\textrm{etc.},\\
\frac{1}{1-n^2z}&=&1+n^2z+n^4z^2+n^6z^3+n^8z^4+\textrm{etc.},\\
\frac{1}{1-n^3z}&=&1+n^3z+n^6z^2+n^9z^3+n^{12}z^4+\textrm{etc.},\\
\frac{1}{1-n^4z}&=&1+n^4z+n^8z^2+n^{12}z^3+n^{16}z^4+\textrm{etc.}\\
&&\textrm{etc.}
\end{eqnarray*}
On the other hand, by putting $-z$ in place of $z$ the series 
$\frac{1}{R}$
follows in a similar way.

30. From the generation of these series it is clear that:
\[
\textrm{I.} \quad \mathfrak{A}=n+n^2+n^3+n^4+n^5+\textrm{etc.},
\]
which is a geometric progression where all powers of $n$ are
multiplied by the coefficient $+1$.

\[
\textrm{II.} \quad \mathfrak{B}=n^2+n^3+2n^4+2n^5+3n^6+3n^7+4n^8+4n^9+\textrm{etc.},
\]
in which the coefficient of each power of $n$ contains as much unities as there
are different ways in which the exponent of $n$ can be partitioned into
two parts, either equal or unequal. Thus the coefficient of the power
$n^8$ is $4$, because $8$ can be partitioned into $2$ parts in four ways
\[
8=1+7, \quad 8=2+6, \quad 8=3+5, \quad 8=4+4.
\]

\[
\textrm{III.} \quad \mathfrak{C}=n^3+n^4+2n^5+3n^6+4n^7+5n^8+7n^9+\textrm{etc.},
\]
in which the coefficient of each power of $n$ contains as many unities
as there are different ways in which the exponent of $n$ can be distributed
into three parts, either equal or unequal. Thus $n^9$ has the coefficient
$7$, because $7$ permits its separation into three parts in $9$ ways:
\[
\begin{array}{ll}
9=1+1+7,&9=1+4+4,\\
9=1+2+6,&9=2+2+5,\\
9=1+3+5,&9=2+3+4,
\end{array}
\]
\[
9=3+3+3.
\]

\[
\textrm{IV.} \quad \mathfrak{D}=n^4+n^5+2n^6+3n^7+5n^8+6n^9+9n^{10}+\textrm{etc.},
\]
where the coefficient of any power of $n$ contains as many unities as
there are different ways in which the exponent of $n$ can be resolved
into four parts, either equal or unequal.
And there is a similar rule for the following series which are found
for the letters $\mathfrak{E},\mathfrak{F},\mathfrak{G}$ etc.

31. The other problem which the Insightful Naud\'e proposed to me along with
the preceding can thus also be resolved by means of these series. It can be
expressed thus:

{\em To find how many different way a given number $m$ can be partitioned
into $\mu$ parts, either equal or unequal, or to find how many different ways
a given number $m$ can be produced by addition from $\mu$ integral
numbers, either equal or unequal.}

The difference between this problem and the preceding is that in the preceding,
the partition was restricted just to parts that were mutually unequal, while
here equal parts are also allowed. 
For expressing the number of all these ways in this problem, the search
will use this form
\[
m^{(\mu)},
\]
which declares namely how many different ways the number $m$ can be partitioned
into $\mu$ integral parts, with equality of some not excluded;
for the letter $i$ which was previously affixed to the above sign $(\mu)$,
which indicated unequal parts, is omitted here.

32. The solution of this problem can be reduced thus to the formation
of the series $\mathfrak{A},\mathfrak{B},\mathfrak{C},\mathfrak{D},
\mathfrak{E}$ etc.; we have already shown above (\S 5) the values
of these letters can be defined from the already known values
of the letters $\alpha,\beta,\gamma,\delta$ etc.
Though this method is general and attacks the very nature of the problem,
it could however seem that the law by which these values proceed is not clear enough
to the eyes. Therefore I shall investigate
the values of the letters $\mathfrak{A},\mathfrak{B},
\mathfrak{C},\mathfrak{D},\mathfrak{E}$ etc.in this case by a means similar
to what I used above (\S 21).

Since it is
\[
\frac{1}{S}=\frac{1}{(1-nz)(1-n^2z)(1-n^3z)(1-n^4z)\, \textrm{etc.}},
\]
it is clear that writing $nz$ in place of $z$ in this form will produce
this form
\[
\frac{1}{(1-n^2z)(1-n^3z)(1-n^4z)(1-n^5z)\, \textrm{etc.}}.
\]
But the prior form $\frac{1}{S}$ is turned into this 
if it is multiplied by $1-nz$. Whence, since we have assumed that
\[
\frac{1}{S}=1+\mathfrak{A}z+\mathfrak{B}z^2+\mathfrak{C}z^3+\mathfrak{D}z^4+\mathfrak{E}z^5+\textrm{etc.},
\]
we put here $nz$ in place of $z$ and we will have
\[
1+\mathfrak{A}nz+\mathfrak{B}n^2z^2+\mathfrak{C}n^3z^3+\mathfrak{D}n^4z^4+\textrm{etc.}
\]
Now let us multiply the first series $\frac{1}{S}$ by $1-nz$
\[
\begin{array}{rrrrrr}
1&+\mathfrak{A}z&+\mathfrak{B}z^2&+\mathfrak{C}z^3&+\mathfrak{D}z^4&+\textrm{etc.}\\
&-nz&-\mathfrak{A}nz^2&-\mathfrak{B}nz^3&-\mathfrak{C}nz^4&-\textrm{etc.}
\end{array}
\]
Since this form should be equal to the previous, it will be
\[
\begin{split}
&\mathfrak{A}=\frac{n}{1-n}=\frac{n}{1-n},\\
&\mathfrak{B}=\frac{\mathfrak{A}n}{1-n^2}=\frac{n^2}{(1-n)(1-n^2)},\\
&\mathfrak{C}=\frac{\mathfrak{B}n}{1-n^3}=\frac{n^3}{(1-n)(1-n^2)(1-n^3)},\\
&\mathfrak{D}=\frac{\mathfrak{C}n}{1-n^4}=\frac{n^4}{(1-n)(1-n^2)(1-n^3)(1-n^4)}\\
&\textrm{etc.}
\end{split}
\]

33. Hence a new relation is perceived here between the values of
the letters $\mathfrak{A},\mathfrak{B},\mathfrak{C},\mathfrak{D}$ etc.
and the values of the letters $\alpha,\beta,\gamma,\delta$ etc. It is noteworthy
that these values do not disagree with each other. For collecting \S 21 it is
understood to be
\begin{eqnarray*}
\alpha&=&\mathfrak{A},\\
\beta&=&n\mathfrak{B},\\
\gamma&=&n^3\mathfrak{C},\\
\delta&=&n^6\mathfrak{D},\\
\epsilon&=&n^{10}\mathfrak{E}\\
\textrm{etc.}&&
\end{eqnarray*}
Thus it is clear from the rule of the coefficients that the series 
$\mathfrak{A},\mathfrak{B},\mathfrak{C},\mathfrak{D}$ etc. agree completely
with the series $\alpha,\beta,\gamma,\delta$ etc. except for
the sole difference 
being in the exponents of $n$.
Indeed in the series $\mathfrak{A}$ the exponents are also equal to the exponents
in the series $\alpha$, but in the series $\mathfrak{B}$ the exponents 
are one less than the exponents in the series $\beta$, in the series
$\mathfrak{C}$ the exponents are three less than the exponents in the series
$\gamma$, and the differences proceed according to the trigonal numbers and so on.

34. Therefore by the series $\alpha,\beta,\gamma,\delta$ etc.
which we showed above how to form and by which the first problem of
Naud\'e was resolved, the latter problem proposed by
Naud\'e
can simultaneously be resolved, so that its solution is reduced to
the solution of the first. Namely it will be
\begin{eqnarray*}
m^{(1)}&=&m^{(1)i},\\
m^{(2)}&=&(m+1)^{(2)i},\\
m^{(3)}&=&(m+3)^{(3)i},\\
m^{(4)}&=&(m+6)^{(4)i}
\end{eqnarray*}
and generally
\[
m^{(\mu)}=\Big(m+\frac{\mu(\mu-1)}{2}\Big)^{(\mu)i}
\]
and in turn
\[
m^{(\mu)i}=\Big(m-\frac{\mu(\mu-1)}{2}\Big)^{(\mu)}.
\]
Furthermore, because we have also found that
\[
(m+\mu)^{(\mu)i}=m^{(\mu)i}+m^{(\mu-1)i},
\]
and reducing this to the present case it will be
\[
\Big(m-\frac{\mu(\mu-3)}{2}\Big)^{(\mu)}=\Big(m-\frac{\mu(\mu-1)}{2}\Big)^{(\mu)}
+\Big(m-\frac{(\mu-1)(\mu-2)}{2}\Big)^{(\mu-1)}
\]
or for convenience
\[
m^{(\mu)}=(m-\mu)^{(\mu)}+(m-1)^{(\mu-1)}.
\]
The series for the letters $\mathfrak{A},\mathfrak{B},\mathfrak{C}$ etc.
are easily formed from this property, and so the latter problem is resolved.

35. For an example of this problem that Insightful Man presented
the question of determining in how many different ways the number $50$ can be 
separated into exactly seven parts, either equal or unequal. This question
can hence be reduced to the first problem, with $m=50$ and $\mu=7$, if we search
for how many different ways the number $50+21$ or $71$ can be separated
into seven unequal parts. This in fact can be done in $8946$ different
ways. Indeed besides this, the same number $8946$ indicates (\S 28) how many
different ways $71-28=43$ can be produced by addition from the numbers $1,2,3,4,5,6,7$.
And generally the number of ways $m^{(\mu)}$ in which the number $m$
is resolved into $\mu$ parts, either equal or unequal, also shows
how many different ways the number $m-\mu$ can be produced by addition from
the particular numbers
\[
1,2,3,4,5,\ldots, \mu.
\]

36. At the end of this paper there is a noteworthy observation to make, which
however I have not yet been able to demonstrate with geometric 
rigor. Namely I have observed that if the infinitely many factors
of the product
\[
(1-n)(1-n^2)(1-n^3)(1-n^4)(1-n^5)\, \textrm{etc.},
\]
are expanded by actual multiplication, they produce this series
\[
1-n-n^2+n^5+n^7-n^{12}-n^{15}+n^{22}+n^{26}-n^{35}-n^{40}+n^{51}+\textrm{etc.},
\]
where only those those powers of $n$ occur whose exponents are contained
in the form $\frac{3xx\pm x}{2}$. And if $x$ is an odd number,
the powers of $n$, which are $n^{\frac{3xx\pm x}{2}}$, will have
the coefficient $-1$, while if $x$ is an even number then
the powers $n^{\frac{3xx \pm x}{2}}$ will have the coefficient $+1$.

37. It is also worth noting that the reciprocal series of this, which arises
from the expansion of this fraction
\[
\frac{1}{(1-n)(1-n^2)(1-n^3)(1-n^4)(1-n^5)\, \textrm{etc.}},
\]
yields namely the recurrent series
\[
1+1n+2n^2+3n^3+5n^4+7n^5+11n^6+15n^7+22n^8+\textrm{etc.}
\]
Of course this series multiplied by the above series
\[
1-n-n^2+n^5+n^7-n^{12}-n^{15}+n^{22}+n^{26}-\textrm{etc.}
\]
produces unity. And in the first series the coefficient of any power
of $n$ contains as many unities as there are different ways in which
the exponent of $n$ can be distributed into parts; thus $5$ can be resolved
in seven ways into parts, as
\[
\begin{array}{lll}
5=5,&5=3+2,&5=2+2+1,\\
5=4+1,&5=3+1+1,&5=2+1+1+1,\\
&5=1+1+1+1+1;&
\end{array}
\]
where namely neither the number of parts or inequalities are prescribed.

\end{document}